\documentclass[final, nomarks]{dmtcs-episciences}

\usepackage[utf8]{inputenc}
\usepackage{subfigure}
\usepackage{amsmath,amsthm,amsfonts,amssymb,amscd}
\usepackage{tikz} 
\usepackage{calrsfs}
\usepackage{mathrsfs}

\newtheorem{theorem}{Theorem}

\author[Helmut Prodinger]{Helmut Prodinger\affiliationmark{1,2}}
\title{The height of skew Dyck paths with two variants of downsteps}
\affiliation{
	Department of Mathematics, University of Stellenbosch, Stellenbosch, South Africa
\\
  	NITheCS (National Institute for
  Theoretical and Computational Sciences), Stellenbosch, South Africa.
}
\begin{document}
\publicationdata{vol. 28:3}{2026}{4}{10.46298/dmtcs.17510}{2026-02-14; 2026-02-14; 2026-07-14}{2026-07-14}

\maketitle
\begin{abstract}
	Recently, in the context of walks of hexagonal circle packings, interest has emerged in the family of skew Dyck paths with two variants
	of down-steps. These paths have steps $U, D_g, D_b, L=D_r$. Using generating functions, the kernel method and (in)finite linear systems, contributions to the (average) height and other enumerations are made. As in many similar instances, the average height is of order $\sqrt n$.
	\end{abstract}

\section{Skew Dyck paths}
Dyck paths consist of up-steps $U=(1,1)$ and down-steps $D=(1,-1)$, start at the origin, never go below the $x$-axis and eventually return to the $x$-axis. 
This concept was and still is very popular, and since such paths are enumerated by Catalan numbers, they figure prominently in \cite{Stanley}.

\emph{Skew} Dyck paths are a variation of the above; they consist of up-steps $U=(1,1)$ and down-steps $D=(1,-1)$, never go below the $x$-axis and return to the $x$-axis. 
Additionally, there is  left step $L=(-1,-1)$, which justifies the epitheton \emph{skew}.  A left step is not allowed to follow an up-step, and an up-step is not allowed to
follow a left step. The following decomposition appears already in the pioneering paper \cite{Emor}:

\begin{figure}[ht]
	\begin{center}
	\begin{tikzpicture}[scale=0.2]
		
		
		\path (-0.3,0.3) node                     {$\mathcal{S}$}
		(2,0) node                     {$=$}
		(4,0) node                     {$\varepsilon$}
		(6,0) node                     {$+$}
		(21,0) node                     {$+$}
		;
		\draw (12.3,1.8) node                     {$\mathcal{S}$};
		\draw (8,-1)-- (10,1);
		\draw (14,1)-- (16,-1);
		\draw (18,0) node                     {$\mathcal{S}$};
		\draw (23,-1)-- (25,1);
		\draw (35,1)-- (33,-1);
		\draw (18,0) node                     {$\mathcal{S}$};
		\draw (29.5,1.8) node                     {$\mathcal{S}\setminus\{\varepsilon\}$};
	\end{tikzpicture}
	
	\caption{Symbolic equation for the family of skew Dyck paths.}
	\end{center}
\end{figure}
From this\footnote{One referee suggests  an alternative notation  $\mathcal{S}=\varepsilon+U\mathcal{S}D\mathcal{S}+U(\mathcal{S}\setminus\{\varepsilon\})L$; both notations are fine.} the generating function $S(z)$, where $z$ marks the number of up-steps, is given by
\begin{equation*}
	S=1+zS^2+z(S-1), \quad\text{or}\quad S(z)=\frac{1-z-\sqrt{1-6z+5z^2}}{2z}.
\end{equation*}

Let us recall that we have up-steps and (possibly) coloured down-steps, denoted by $U, D_g, D_b, L=D_r$.

Now we move to an extended family of skew Dyck paths, where two types of down-steps can occur.
They will be marked in green resp.\ blue.
\begin{figure}[ht]
	\begin{center}
	\begin{tikzpicture}[scale=0.2]
		
		
		\path (-0.3,0.3) node                     {$\mathcal{S}$}
		(2,0) node                     {$=$}
		(4,0) node                     {$\varepsilon$}
		(6,0) node                     {$+$}
		(21,0) node                     {$+$}
		;
		\draw (12.3,1.8) node                     {$\mathcal{S}$};
		\draw (8,-1)-- (10,1);
		\draw[thick,green] (14,1)-- (16,-1);
		\draw (18,0) node                     {$\mathcal{S}$};
		
		\draw[xshift=15cm] (12.3,1.8) node                     {$\mathcal{S}$};
		\draw [xshift=15cm](8,-1)-- (10,1);
		\draw[xshift=15cm,thick,blue] (14,1)-- (16,-1);
		\draw [xshift=15cm](18,0) node                     {$\mathcal{S}$};

		\draw[xshift=16cm] (23,-1)-- (25,1);
		\draw [xshift=16cm](35,1)-- (33,-1);
		\draw [xshift=15cm](21.5,0) node                     {$+$};
		\draw [xshift=16cm](29.5,1.8) node                     {$\mathcal{S}\setminus\{\varepsilon\}$};

	\end{tikzpicture}
	
	\caption{Symbolic equation for the family of skew Dyck paths with two variants of down-steps.}
	\end{center}
\end{figure}
The generating function is then
\begin{equation*}
	S=1+\mathbf{2}zS^2+z(S-1), \quad\text{or}\quad S(z)=\frac{1-z-\sqrt{1-10z+9z^2}}{4z}.
\end{equation*}
The sequence of coefficients of $S(z)$  is sequence A086871 in \cite{OEIS}. 

Note that the dominant singularity of this is at $z=\varrho=\frac19$ which leads to an exponential growth order $9^n$; a more precise statement follows later.

My interest in this particular family of skew Dyck paths was triggered by the recent paper \cite{Baril}, where it was demonstrated that
they are in bijection with  a family of walks in hexagonal circle packings. 

Here, we will consider the height of such skew Dyck paths, together with some additional observations. The height is maximal level that the path reaches, which is an integer $\ge0$. We will however, not use the language of walks in hexagonal circle packings.

A general reference about paths in lattice structures is \cite{BF}; various recent variations appear in \cite{garden}.

\section{Enumeration of skew Dyck paths with two variants of down-steps}

The generating function 
\begin{equation*}
	S(z)=\frac{1-z-\sqrt{1-10z+9z^2}}{4z}
\end{equation*}
becomes nicer using the substitution $z=\dfrac{u}{(1+u)(1+4u)}$, since it is simply $S=1+2u$. Apparently, this handy substitution wasn't observed before. To read off coefficients, we use Cauchy's integral formula, see \cite{FO}:
\begin{align*}
	[z^n]S&=\frac1{2\pi i}\oint \frac{dz}{z^{n+1}}S\\&=\frac1{2\pi i}\oint \frac{du(1-4u^2)}{(1+u)^2(1+4u)^2}\frac{(1+u)^{n+1}(1+4u)^{n+1}}{u^{n+1}}(1+2u)\\
	&=\frac1{2\pi i}\oint du(1-4u^2)(1+2u)\frac{(1+u)^{n-1}(1+4u)^{n-1}}{u^{n+1}}\\
	&=[u^n](1+2u-4u^2-8u^3)(1+u)^{n-1}(1+4u)^{n-1}\\
	&=\binom{n-1;1,5,4}{n}+2\binom{n-1;1,5,4}{n-1}-4\binom{n-1;1,5,4}{n-2}-8\binom{n-1;1,5,4}{n-3}.
\end{align*}
We used the notion of weighted trinomial coefficients
\begin{equation*}
	\binom{m;1,5,4}{k}=[v^k](1+5v+4v^2)^m=[v^k](1+v)^m(1+4v)^m=\sum_{j=0}^k4^j\binom{m}{j}\binom{m}{k-j}.
\end{equation*}
The notion of (weighted) trinomial coefficients has successfully appeared in \cite{HPW, prodinger-three} as well as in some papers by Panny and coauthors, see for instance
\cite{panny}.

We found it more intuitive to draw a left step $(-1,-1)$ as a marked (in red) down-step $(1,-1)$, see \cite{korea}. In this setting, we have now Dyck paths with 3 types of down-steps (blue, green, red), say $D_b, D_g, D_r$, with the restriction that $UD_r$ and $D_rU$ are not allowed to occur.

\begin{figure}[ht]
	\begin{center}
	\begin{tikzpicture}[scale=0.3]
		\draw[help lines] (0,0) grid (29,8);
		
		\draw [thick] (0,0) -- (1,1)--(2,2);
		
		\draw [thick,green] (2,2)--(3,1);
		\draw [thick,blue] (3,1)--(4,0);
		\draw [thick] (4,0)--(10,6);
		\draw [thick,blue] (10,6)--(11,5);
		\draw [thick,red] (11,5)--(13,3);
		\draw [thick,green] (13,3)--(14,2);
		\draw [thick] (14,2)--(19,7);
		\draw [thick,green] (19,7)--(21,5);
		\draw [thick] (21,5)--(22,6);
		\draw [thick,blue] (22,6)--(24,4);
		\draw [thick,green] (24,4)--(27,1);
		\draw [thick,red] (27,1)--(28,0);
		
	\end{tikzpicture}
	
	\caption{A path as studied in the paper of length 28 (=semi-length 14) and height 7.}
	\end{center}
\end{figure}

We will  construct an automaton that checks whether the syntactic conditions are still fulfilled, i. e., $UL$ and $LU$ are forbidden. 
According to the three types of steps, we have three layers of states (current level when scanning the word from left to right).

It will become important soon that the dominant singularity $\varrho=\frac19$ corresponds to the value $u_{\varrho}=\frac{1}{2}$.

The decomposition for the family $\mathcal{S}$ translates as well to $\mathcal{S}_h$ of bounded height, viz.
\begin{equation*}
	S_h=1+2zS_{h-1}S_h+z(S_{h-1}-1)
\end{equation*}
or
\begin{equation*}
	S_h(1-2zS_{h-1})=1-z+zS_{h-1},
\end{equation*}
resp.
\begin{equation*}
	S_h=\frac{1-z+zS_{h-1}}{1-2zS_{h-1}}.
\end{equation*}
This recursion will be solved using the \emph{ansatz} $S_h=\frac{p_h}{q_h}$:
\begin{equation*}
	\frac{p_h}{q_h}=\frac{1-z+z\frac{p_{h-1}}{q_{h-1}}}{1-2z\frac{p_{h-1}}{q_{h-1}}}=
	\frac{(1-z)q_{h-1}+zp_{h-1}}{q_{h-1}-2zp_{h-1}}.
\end{equation*}
Comparing numerators and denominators we derive
\begin{align*}
	p_h&=(1-z)q_{h-1}+zp_{h-1},\\
	q_h&=q_{h-1}-2zp_{h-1};
\end{align*}
this system may be solved using the previous substitution
\begin{gather*}
	p_h=\frac{-2u}{1-2u}\bigg(\frac{(4u+3)u}{(1+u)(1+4u)}\bigg)^h+\frac1{1-2u}\bigg(\frac{3u+1}{(1+u)(1+4u)}\bigg)^h,\\
	q_h=\frac{-4u^2}{1-4u^2}\bigg(\frac{(4u+3)u}{(1+u)(1+4u)}\bigg)^h+\frac1{1-4u^2}\bigg(\frac{3u+1}{(1+u)(1+4u)}\bigg)^h.
\end{gather*}
Although a computer was used to get this, one could check that these answers indeed solve the recursions. Now
\begin{equation*}
	\frac{p_h}{q_h}=(1+2u)\frac{-2u(4u+3)^hu^h+(3u+1)^h}{-4u^2(4u+3)^hu^h+(3u+1)^h}.
\end{equation*}
In the limit $h\to\infty$ (no height restrictions) $S_\infty=\lim_{h\to\infty}S_h=1+2u$, as it should. 
Then we set $S^{[>h]}:=S_\infty-S_h$ and
\begin{align*}
	S^{[>h]}&:=1+2u-\frac{p_h}{q_h}=(1+2u)\bigg[1-\frac{-2u(4u+3)^hu^h+(3u+1)^h}{-4u^2(4u+3)^hu^h+(3u+1)^h}\bigg]\\
	&=(1+2u)\frac{-4u^2(4u+3)^hu^h+(3u+1)^h+2u(4u+3)^hu^h-(3u+1)^h}{-4u^2(4u+3)^hu^h+(3u+1)^h}\\
	&=(1-4u^2)2u\frac{(4u+3)^hu^h}{-4u^2(4u+3)^hu^h+(3u+1)^h}\\
	&=2u(1-4u^2)\frac{\lambda^h}{1-4u^2\lambda^h}=\frac{1-4u^2}{2u}\frac{4u^2\lambda^h}{1-4u^2\lambda^h},
\end{align*}
with
\begin{equation*}
	\lambda:=\frac{(4u+3)u}{3u+1}.
\end{equation*}
The critical value $u_{\varrho}=\frac12$ corresponds to $\lambda_{\varrho}=1$. We approximate both, $2u$ and $\lambda$ by $e^{-t}$. Since they are not the same, this will lead to a leading term in the asymptotic expansion only.
In order to get the average height of our family of skew Dyck paths we need to study
\begin{equation*}
	\sum_{h\ge0}S^{[>h]}=\frac{1-4u^2}{2u}\sum_{h\ge1}\frac{4u^2\lambda^h}{1-4u^2\lambda^h}
\end{equation*}
or, approximating it,
\begin{equation*}
	\frac{1-4u^2}{2u}\sum_{h\ge1}\frac{e^{-(h+2)t}}{1-e^{-(h+2)t}}
	\sim(2t+\dots)\sum_{h\ge1}\frac{e^{-ht}}{1-e^{-ht}}.
\end{equation*}
The asymptotic expansion of the series appeared in many papers, in particular in  \cite{HPW}, where many technical details are provided. We continue with
\begin{equation*}
	\sim\frac{1-\lambda^2}{\lambda}\frac{-\log(1-\lambda)}{1-\lambda}\sim -\frac{1+\lambda}{\lambda}\log(1-\lambda).
\end{equation*}
We need the expansion around the critical values;
\begin{equation*}
	\sqrt{1-9z}\sim \frac{\sqrt2}{3}(1-2u),
\end{equation*}
\begin{equation*}
	1-\lambda\sim\frac45(1-2u)\sim \frac45\frac{3}{\sqrt2}\sqrt{1-9z}=\frac{6\sqrt2}{5}\sqrt{1-9z},
\end{equation*}
\begin{equation*}
	\log(1-\lambda)\sim\log\bigg(\frac{6\sqrt2}{5}\sqrt{1-9z}\bigg)\sim\frac12\log(1-9z),
\end{equation*}
\begin{equation*}
	-\frac{1+\lambda}{\lambda}\log(1-\lambda)\sim (-2) \frac12\log(1-9z)=-\log(1-9z).
\end{equation*}
By singularity analysis, the coefficient of $z^n$ in the last approximation is asymptotic to $\dfrac{9^n}{n}$. For the average height, we need to divide this by the total number of relevant paths of semi-length $n$.

The generating function $S(z)$ will be expanded:
\begin{equation*}
	S(z)\sim2-\frac{3}{\sqrt2}\sqrt{1-9z}.
\end{equation*}
Thus (singularity analysis, \cite{FO})
\begin{equation*}
	[z^n]S(z)\sim -\frac{3}{\sqrt2}9^n\frac{n^{-3/2}}{\Gamma(-\frac12)}=\frac{3}{2\sqrt2}9^n\frac{n^{-3/2}}{\sqrt{\pi}};
\end{equation*}
the desired quotient is then asymptotic to
\begin{equation*}
	\frac{\frac{9^n}{n}}{\frac{3}{2\sqrt2}9^n\frac{n^{-3/2}}{\sqrt{\pi}}}=\frac{2\sqrt2}{3}\sqrt{\pi n}\sim0.9428090414\sqrt{\pi n}.
\end{equation*}

\begin{theorem}
	The average height of skew Dyck paths with two types of down-steps is \linebreak asymptotic to
	\begin{equation*}
		\frac{2\sqrt2}{3}\sqrt{\pi n}
	\end{equation*}
	where all such paths of semi-length $n$ are equally likely.
	
	This may be compared to the formula of standard skew Dyck paths (see \cite{garden})
	\begin{equation*}
		\frac{2}{\sqrt5}\sqrt{\pi n}
	\end{equation*}
	and for standard Dyck paths (not skew)
	\begin{equation*}
		\sqrt{\pi n}.
	\end{equation*}
	The last result is classical \cite{deBrKnRi}.
\end{theorem}

\section{The  recursions}

We start with an automaton that scans an input word of the alphabet \{up, blue down, green down, red down\} whether is satisfies the necessary syntactical conditions.
It has three layers, according to the type of the last step. Further, we will work with generating functions $f_i$, $g_i$, $h_i$, for $i\ge0$, according to the state in which the path has landed.
Eventually, $f_0+g_0+h_0$ will be the generating function $S(z^2)$; note that for the present approach we need to consider the length, not the semi-length. We will be able
to compute $f_i$, $g_i$, $h_i$ explicitly.

\begin{figure}[ht]
	\begin{center}
	\begin{tikzpicture}[line width=1pt,scale=0.3]
		
		\draw (0,4) circle (0.3cm);
		\fill (0,4) circle (0.3cm);
		
		\foreach \t in{0,1,2,3,4,5,6,7}
		{
			\foreach \x in {0,2}
			{
				\def \y{\x-0};
				\draw[blue](\x+4*\t,\y) -- (\x+1+4*\t,\y+1);
			}
			\foreach \x in {1,3}
			{
				\def \y{\x-0};
				\draw[green](\x+4*\t,\y) -- (\x+1+4*\t,\y+1);
			}
			\draw[blue,latex-](0+4*\t,0) -- (1+4*\t,1);
		}
		
		\foreach \t in{0,1,2,3,4,5,6,7}
		{
			\foreach \x in {0,2}
			{
				\def \y{\x-0};
				\draw[blue](\x+4*\t,0) -- (\x+1+4*\t,0);
			}
			\foreach \x in {1,3}
			{
				\def \y{\x-0};
				\draw[green](\x+4*\t,0) -- (\x+1+4*\t,0);
			}
			\draw[blue,latex-](0+4*\t,0) -- (1+4*\t,0);
		}
		
		\foreach \t in{0,1,2,3,4,5,6,7}
		{\draw[-latex](4*\t,4)--(4*\t+4,4);}
		\foreach \t in{0,1,2,3,4,5,6,7}
		{\draw[-latex](4*\t,0)to [out=90,in=-155](4*\t+4,4);}
		\foreach \t in{0,1,2,3,4,5,6,7}
		{\draw[red,latex-](4*\t,-4)to [out=65,in=-155](4*\t+4,0);\draw[red,latex-](4*\t,-4)to (4*\t+4,-4);}
		
		\foreach \t in{1,2,3,4,5,6,7,8}
		{
			\foreach \x in {0}
			{
				\def \y{\x-3};
				\draw[blue](\x+4*\t-1,\y) -- (\x+4*\t,\y-1);
			}
			\foreach \x in {2}
			{
				\def \y{\x-3};
				\draw[blue](\x+4*\t-1-4,\y) -- (\x+4*\t-4,\y-1);
			}
			\foreach \x in {1}
			{
				\def \y{\x-3};
				\draw[green](\x+4*\t+1-4,\y) -- (\x+4*\t+2-4,\y-1);
			}
			\foreach \x in {3}
			{
				\def \y{\x-3};
				\draw[green, latex-](\x+4*\t+1-8,\y) -- (\x+4*\t+2-8,\y-1);
			}
		}	
		
	\end{tikzpicture}
	\caption{An automaton to check the syntax of skew Dyck paths with two types of down-steps (blue resp. green); the left step is drawn in red. }
	\label{pso}
	\end{center}
\end{figure}

The system of recursions can be directly read off, by considering the last step:
\begin{align*}
	f_0&=1,\quad f_{i+1}=zf_i+zg_i,\\
	g_i&=2zf_{i+1}+2zg_{i+1}+2zh_{i+1},\\
	h_i&=zg_{i+1}+zh_{i+1}.
\end{align*}
The first method of solving this is using bivariate generating functions, viz.
\begin{equation*}
	F(u)=F(u,z)=\sum_{i\ge0}u^if_i,\
	G(u)=G(u,z)=\sum_{i\ge0}u^ig_i,\
	H(u)=H(u,z)=\sum_{i\ge0}u^ih_i.
\end{equation*}
From the recursions,
\begin{align*}
	\sum_{i\ge0}u^{i+1}f_{i+1}&=\sum_{i\ge0}u^{i+1}zf_i+\sum_{i\ge0}u^{i+1}zg_i,\\*
	\sum_{i\ge0}u^{i+1}g_i&=\sum_{i\ge0}u^{i+1}2zf_{i+1}+\sum_{i\ge0}u^{i+1}2zg_{i+1}+\sum_{i\ge0}u^{i+1}2zh_{i+1},\\*
	\sum_{i\ge0}u^{i+1}h_i&=\sum_{i\ge0}u^{i+1}zg_{i+1}+\sum_{i\ge0}u^{i+1}zh_{i+1},
\end{align*}
or
\begin{align*}
	F(u)-1&=zuF(u)+zuG(u),\\
	uG(u)&=2z(F(u)-1)+2z(G(u)-G(0))+2z(H(u)-H(0)),\\
	uH(u)&=z(G(u)-G(0))+z(H(u)-H(0)).
\end{align*}
This system can be solved:
\begin{align*}
	F(u)&=\frac{-z(-3+2z^2)+(2z^2+2z^2G(0)+2z^2H(0)-1)u}{zu^2-(1+z^2)u+z(3-z^2)},\\
	G(u)&=\frac{z(G(0)+H(0)+z^2)-2z^2(G(0)+H(0)+1)u}{zu^2-(1+z^2)u+z(3-z^2)},\\
	H(u)&=\frac{-z(2z^2H(0)+2z^2+2z^2G(0)-H(0)-G(0)-z^2(H(0)+G(0))u}{zu^2-(1+z^2)u+z(3-z^2)}.
\end{align*}
The denominator will be factored:
\begin{equation*}
	zu^2-(1+z^2)u+z(3-z^2)=z(u-\mu_1)(u-\mu_2)
\end{equation*}
with
\begin{align*}
	\mu_1&=\frac{1+z^2+\sqrt{1-10z^2+9z^4}}{2z},\\
	\mu_2&=\frac{1+z^2-\sqrt{1-10z^2+9z^4}}{2z}.
\end{align*}
In the spirit of the kernel method \cite{garden} the factor $u-\mu_2$ may be divided out from numerators and denominators, with the result
\begin{align*}
	F(u)&=\frac{-2{z}^{2}g_0-2{z}^{2}h_0-2{z}^{2}+1}{z(u-\mu_1)},\\
	G(u)&=\frac{2z ( g_0+h_0+1) }{z(u-\mu_1)},\\
	H(u)&=\frac{z ( g_0+h_0) }{z(u-\mu_1)}.
\end{align*}
Now one can plug in $u=0$ and identify $f_0=1$, $g_0$, $h_0$;
\begin{align*}
	G(0)=g_0&=\frac{2z(-z+\mu_1)}{\mu_1(-3z+\mu_1)},\\
	H(0)=h_0&=\frac{2z^2}{\mu_1(-3z+\mu_1)}.
\end{align*}
Hence
\begin{align*}
	F(u)&=\frac{-2{z}^{2}g_0-2{z}^{2}h_0-2{z}^{2}+1}{z(u-\mu_1)}=\frac{-\mu_1}{u-\mu_1}=\frac{1}{1-u/\mu_1},\\
	G(u)&=\frac{2z ( g_0+h_0+1) }{z(u-\mu_1)}=\frac{\mu_1-\frac1z}{u-\mu_1}=\frac{\frac1{z\mu_1}-1}{1-u/\mu_1},\\
	H(u)&=\frac{z ( g_0+h_0) }{z(u-\mu_1)}=\frac{\frac{\mu_1}2+z-\frac1{2z}}{u-\mu_1}=\frac{-\frac{1}2-(z-\frac1{2z})\frac1{\mu_1}}{1-u/\mu_1}
\end{align*}
and further
\begin{align*}
	[u^k]F(u)&=\frac1{\mu_1^k},\\
	[u^k]G(u)&=\frac1{z\mu_1^{k+1}}-\frac1{\mu_1^{k}}\\
	[u^k]H(u)&=-\frac{1}2\frac1{\mu_1^k}-\Big(z-\frac1{2z}\Big)\frac1{\mu_1^{k+1}}.
\end{align*}
\begin{theorem}
	The generating functions of skew Dyck paths (or prefixes of them) with two variants of down-steps leading to level $i$ are given by
	\begin{align*}
		f_i&=\frac1{\mu_1^i},\\*
		g_i&=\frac1{z\mu_1^{i+1}}-\frac1{\mu_1^{i}},\\*
		h_i&=-\frac{1}2\frac1{\mu_1^i}-\Big(z-\frac1{2z}\Big)\frac1{\mu_1^{i+1}}.
	\end{align*}
	These may be rewritten using $\mu_2=\dfrac{3-2z^2}{\mu_1}$. The total generating function leading to level $i$ is
	\begin{equation*}
		f_i+g_i+h_i=\frac1{\mu_1^i}+\frac1{z\mu_1^{i+1}}-\frac1{\mu_1^{i}}-\frac{1}2\frac1{\mu_1^i}-\Big(z-\frac1{2z}\Big)\frac1{\mu_1^{i+1}}=
		-\frac{1}2\frac1{\mu_1^i}-\Big(z-\frac3{2z}\Big)\frac1{\mu_1^{i+1}}.
	\end{equation*}
	In particular
	\begin{equation*}
		f_0+g_0+h_0=
		\frac{1-z^2-\sqrt{1-9z^2+10z^4}}{4z^2}.
	\end{equation*}
\end{theorem}
The quantity $F(1)+G(1)+H(1)$ is also of interest since it is the generating function of \emph{any} prefix of our class of skew Dyck paths. The series is
\begin{equation*}
	F(1)+G(1)+H(1)=1+z+3{z}^{2}+5{z}^{3}+17{z}^{4}+31{z}^{5}+107{z}^{6}+209{z}^{7}+725{z}^{8}+1483{z}^{9}+5159{z}^{10}+\dots
\end{equation*}
which is not in \cite{OEIS}.

The denominator $zu^2-(1+z^2)u+z(3-z^2)$ drives the recursions ($i\ge0$)
\begin{align*}
	zf_i-(1+z^2)f_{i+1}+z(3-2z^2)f_{i+2}&=0,\\
	zg_i-(1+z^2)g_{i+1}+z(3-2z^2)g_{i+2}&=0,\\
	zh_i-(1+z^2)h_{i+1}+z(3-2z^2)h_{i+2}&=0;
\end{align*}
the quantities $f_0,g_0, h_0$ are known, and further
\begin{align*}
	-(1+z^2)f_{0}+z(3-z^2)f_{1}&=-z\mu_1,\\
	(2-3z^2)g_0+z(3-2z^2)g_1&=\frac{2(1-z^2)}{z},\\
	-(2z+1)(2z-1)h_0-z(2z^2-1)(2z^2-3)h_1&=2z^4.
\end{align*}
The last relations may be checked directly.
Now we can express the sequences using matrices.
\begin{equation*}
	\begin{pmatrix}
		-(1+z^2) & z(3-2z^2)& 0&\cdots    \\
		z & -(1+z^2) &z(3-2z^2)& \cdots   \\
		&		z & -(1+z^2) &z(3-2z^2)   \\
		\vdots  & \ddots  & \ddots & \ddots   
	\end{pmatrix}
	\begin{pmatrix}
		f_0\\
		f_1\\
		f_2\\
		\vdots
	\end{pmatrix}=
	\begin{pmatrix}
		-z\mu_1\\
		0\\
		0\\
		\vdots
	\end{pmatrix}
\end{equation*}

\begin{equation*}
	\begin{pmatrix}
		(2-3z^2) & z(3-2z^2)& 0&\cdots    \\
		z & -(1+z^2) &z(3-2z^2)& \cdots   \\
		&		z & -(1+z^2) &z(3-2z^2)   \\
		\vdots  & \ddots  & \ddots & \ddots  
	\end{pmatrix}
	\begin{pmatrix}
		g_0\\
		g_1\\
		g_2\\
		\vdots
	\end{pmatrix}=
	\begin{pmatrix}
		\tfrac{2(1-z^2)}{z}\\
		0\\
		0\\
		\vdots
	\end{pmatrix}
\end{equation*}

\begin{equation*}
	\begin{pmatrix}
		-(1-4z^2) & z(3-2z^2)(1-2z^2)& 0&\cdots    \\
		z & -(1+z^2) &z(3-2z^2)& \cdots   \\
		&		z & -(1+z^2) &z(3-2z^2)   \\
		\vdots  & \ddots  & \ddots & \ddots  
	\end{pmatrix}
	\begin{pmatrix}
		h_0\\
		h_1\\
		h_2\\
		\vdots
	\end{pmatrix}=
	\begin{pmatrix}
		-2z^4\\
		0\\
		0\\
		\vdots
	\end{pmatrix}.
\end{equation*}
Now we will only work out $f_0, g_0, h_0$, provided that the matrix is cut off after $K+1$ rows/columns, i.e. $f_i:=0$ for $i>K$ and similarly $g_i:=0$ for $i>K$ and
$h_i:=0$ for $i>K$. Since $f_0=1$ independently of $K$ we only consider the $g$-matrix and the $h$-matrix.
Let us use the notations $f_0^{[K]}, g_0^{[K]}, h_0^{[K]}$. Since $f_0^{[K]}=1$, for any $K$, we don't need to compute it. For the remaining two quantities, we solve a linear system using Cramer's rule. So we have to deal with various determinants. Consider
\begin{equation*}
	\begin{pmatrix}
		-(1+z^2) & z(3-2z^2)& 0&\cdots  &0  \\
		z & -(1+z^2) &z(3-2z^2)& \cdots & 0 \\
		&		z & -(1+z^2) &z(3-2z^2) & 0 \\
		\vdots  & \ddots  & \ddots & \ddots  & \vdots\\
		&		&& & -(1+z^2)  \\
	\end{pmatrix}
\end{equation*}
with  $K$ rows and $\mathcal{D}_K$ its  determinant. Then
$\mathcal{D}_0=1$, $\mathcal{D}_1=-(1+z^2)$, and by expansion
\begin{equation*}
	\mathcal{D}_K=-(1+z^2)\mathcal{D}_{K-1}-z^2(3-2z^2)\mathcal{D}_{K-2}.
\end{equation*}
This can be solved by standard methods;
\begin{equation*}
	\mathcal{D}_K=\frac{1}{\sqrt{1-10z^2+9z^4}}\Big[(-z\mu_2)^{K+1}-(-z\mu_1)^{K+1}\Big]=
	\frac{(-z)^{K+1}}{\sqrt{1-10z^2+9z^4}}\Big[\mu_2^{K+1}-\mu_1^{K+1}\Big].
\end{equation*}
Notice that
\begin{equation*}
	\frac{\mathcal{D}_{K-1}}{\mathcal{D}_K}=\frac{\frac{(-z)^{K}}{\sqrt{1-10z^2+9z^4}}\Big[\mu_2^{K}-\mu_1^{K}\Big]}{\frac{(-z)^{K+1}}{\sqrt{1-10z^2+9z^4}}\Big[\mu_2^{K+1}-\mu_1^{K+1}\Big]}=-\frac1z\frac{\mu_2^{K}-\mu_1^{K}}{\mu_2^{K+1}-\mu_1^{K+1}}.
\end{equation*}
We consider the truncated system
\begin{equation*}
	\begin{pmatrix}
		(2-3z^2) & z(3-2z^2)& 0&\cdots  &0  \\
		z & -(1+z^2) &z(3-2z^2)& \cdots & 0 \\
		&		z & -(1+z^2) &z(3-2z^2) & 0 \\
		\vdots  & \ddots  & \ddots & \ddots  & \vdots\\
		&&&		z & -(1+z^2)  \\
	\end{pmatrix}
	\begin{pmatrix}
		g_0\\
		g_1\\
		g_2\\
		\vdots\\
		g_K
	\end{pmatrix}=
	\begin{pmatrix}
		\tfrac{2(1-z^2)}{z}\\
		0\\
		0\\
		\vdots\\
		0
	\end{pmatrix}
\end{equation*}
Then
\begin{align*}
	g_0^{[K]}&=\frac{\tfrac{2(1-z^2)}{z}\mathcal{D}_{K}}{(2-3z^2)\mathcal{D}_{K}-z^2(3-2z^2)\mathcal{D}_{K-1}}\\
	&=\frac{\tfrac{2(1-z^2)}{z(2-3z^2)}}{1-\frac{z^2(3-2z^2)}{(2-3z^2)}\frac{\mathcal{D}_{K-1}}{\mathcal{D}_{K}}}.
\end{align*}
From
\begin{equation*}
	\frac{\mathcal{D}_{K-1}}{\mathcal{D}_{K}}=\frac1{-(1+z^2)-z^2(3-2z^2)\frac{\mathcal{D}_{K-2}}{\mathcal{D}_{K-1}}}
\end{equation*}
one could write a continued fraction expansion for $g_0^{[K]}$, see \cite{Shallit-paper} for similar expansions.

One could as well engage into asymptotics as we did earlier with $S^{[h]}$, but we won't do that in order to save space;
the paper \cite{cornerless} has such expansions for the interested readers.

Let as briefly write the analogous developments for $h_0^{[K]}$
\begin{equation*}
	\begin{pmatrix}
		-(1-4z^2) & z(3-2z^2)(1-2z^2)& 0&\cdots  &0  \\
		z & -(1+z^2) &z(3-2z^2)& \cdots & 0 \\
		&		z & -(1+z^2) &z(3-2z^2) & 0 \\
		\vdots  & \ddots  & \ddots & \ddots  & \vdots\\
		&&&		z & -(1+z^2)  \\
	\end{pmatrix}
	\begin{pmatrix}
		h_0\\
		h_1\\
		h_2\\
		\vdots\\
		h_K
	\end{pmatrix}=
	\begin{pmatrix}
		-2z^4\\
		0\\
		0\\
		\vdots\\
		0
	\end{pmatrix};
\end{equation*}

\begin{align*}
	h_0^{[K]}&=\frac{-2z^4\mathcal{D}_{K}}{-(1-4z^2)\mathcal{D}_{K}-z^2(3-2z^2)(1-2z^2)\mathcal{D}_{K-1}}.
\end{align*}

\section{Further research}

Apart from possible considerations about restricted walks in the graph in Figure~\ref{pso} we mention the possibility to allow more than two colors.
Some brief computations seem to indicate that results about this scenario would be less pretty.

\clearpage


\bibliographystyle{plain}

\label{sec:biblio}

\end{document}